\documentclass[reqno]{mfatshort}
\usepackage[english]{babel}
\usepackage{amssymb, latexsym, euscript}
 \newtheorem{theorem}{\bfseries{Theorem}}[section]
 \newtheorem{lemma}{\bfseries{Lemma}}[section]
 \newtheorem{rmk}{\bfseries{Remark}}[section]
 \newtheorem*{rem*}{\bfseries{Remark}}
 
 \newtheorem*{prop*}{\bfseries{Proposition}}
 \newtheorem{cor}{\bfseries{Corollary}}[section]

\newcommand{\HH}{{\mathfrak H}}

\newcommand{\R}{{\mathbb R}}

\newcommand{\limind}{\mathop{{\rm lim\,ind}}\limits}

\begin{document}

\title[DIRECT AND INVERSE THEOREMS ...]{DIRECT AND INVERSE THEOREMS IN THE THEORY
OF APPROXIMATION BY THE RITZ METHOD}
\author{S. M. Torba}
\email{sergiy.torba@gmail.com}
\author{M. L. Gorbachuk}
\email{imath@horbach.kiev.ua}
\author{Ya. I. Grushka}
\email{grushka@imath.kiev.ua}
\address{Institute of Mathematics, Ukrainian Academy of Sciences, Kyiv.}

\subjclass[2000]{Primary 41A25, 41A17, 41A65} \keywords{ Direct and
inverse theorems, modulo of continuity, Hilbert space, entire
vectors of exponential type}

\begin{abstract}
For an arbitrary self-adjoint operator $B$ in a Hilbert space $\HH$,
we present direct and inverse theorems establishing the relationship
between the degree of smoothness of a vector $x\in\HH$ with respect
to the operator $B$, the rate of convergence to zero of its best
approximation by exponential-type entire vectors of the operator
$B$, and the $k$-modulus of continuity of the vector $x$ with
respect to the operator $B$. The results are used for finding a
priori estimates for the Ritz approximate solutions of operator
equations in a Hilbert space.
\end{abstract}

\maketitle

 \section{Introduction}
Let $B$ be a closed linear operator with dense domain of definition
$\mathcal{D}(B)$ in a separable Hilbert space $\mathfrak{H}$ over
the field of complex numbers.

Let $C^{\infty }(B)$ denote the set of all infinitely differentiable
vectors of the operator $B$, i.e.,
\[
C^{\infty }(B)=\bigcap _{n\in
\mathbb{N}_{0}}\mathcal{D}(B^{n}),\quad \mathbb{N}_{0}=\left\{
0,1,2,...\right\} =\mathbb{N}\cup \{0\}.
\]
For a number $\alpha >0$, we set
\[
\mathfrak{E}^{\alpha }(B)=\left\{
x\in C^{\infty }(B)\, |\, \exists c=c(x)>0\, \, \forall k\in
\mathbb{N}_{0}\, \left\Vert B^{k}x\right\Vert \leq c\alpha
^{k}\right\} .
\]

The set $\mathfrak{E}^{\alpha }(B)$ is a Banach space with respect
to the norm
\[
\left\Vert x\right\Vert _{\mathfrak{E}^{\alpha
}(B)}=\sup _{n\in \mathbb{N}_{0}}\frac{\left\Vert B^{n}x\right\Vert
}{\alpha ^{n}}\, .
\]
 Then $\mathfrak{E}(B)=\bigcup \limits _{\alpha >0}\mathfrak{E}^{\alpha }(B)$
is a linear locally convex space with respect to the topology of the
inductive limit of the Banach spaces $\mathfrak{E}^{\alpha }(B)$:
\[
\mathfrak{E}(B)=\limind _{\alpha \rightarrow \infty
}\mathfrak{E}^{\alpha }(B).
\]
 Elements of the space $\mathfrak{E}(B)$ are called exponential-type entire vectors of the operator $B$.~
  The type $\sigma (x,B)$
of a vector $x\in \mathfrak{E}(B)$ is defined as the number
\[
\sigma (x,B)=\inf \left\{ \alpha >0\, :\, x\in \mathfrak{E}^{\alpha
}(B)\right\} =\limsup _{n\rightarrow \infty }\left\Vert
B^{n}x\right\Vert ^{\frac{1}{n}}.
\]
In what follows, we always assume that the operator $B$ is
self-adjoint in $\mathfrak{H}$, and $E(\Delta )$ is its spectral
measure.

Let $G(\cdot )$ be an almost everywhere finite measurable function
on $\R$. A function $G(B)$ of the operator $B$ is understood as
follows:
\[
G(B):=\int _{-\infty }^{\infty }G(\lambda )dE(\lambda ).
\]
As shown in \cite{MGorb_AnSolutions}, one has $\mathfrak{E}^{\alpha
}(B)=E([-\alpha ,\alpha ])\mathfrak{H}$ for every $\alpha >0$.

According to \cite{Kyptsov}, we set
\begin{equation}
\omega
_{k}(t,x,B)=\sup _{0<\tau \leq t}\left\Vert \Delta _{\tau
}^{k}x\right\Vert ,\, \, k\in
\mathbb{N},\label{w_{d}ef}
\end{equation}
 where
\begin{align}
\Delta _{h}^{k}=(U(h)-\mathbb{I})^{k}=\sum
_{j=0}^{k}(-1)^{k-j}C_{k}^{j}U(jh),\quad k\in \mathbb{N}_{0},\, \,
h\in \mathbb{R} \label{Delta_def}\quad(\Delta _{h}^{0}\equiv 1,\, \,
h\in \mathbb{R}_+),
\end{align}
 and $U(h)=\exp (ihB)$ is the group of unitary operators in $\mathfrak{H}$
with generator $iB$ \cite{Axizer-Glazman}.

The definition of $\omega _{k}(t,x,B)$ implies that the following
assertions are true $k\in \mathbb{N}$:

\begin{enumerate}
\item $\omega _{k}(0,x,B)=0$;
\item for fixed $x$, the function $\omega _{k}(t,x,B)$ does not decrease on $\mathbb{R}_{+}=[0,\infty )$;
\item $\omega _{k}(\alpha t,x,B)\leq [1+\alpha ]^{k}\omega _{k}(t,x,B)$
~ ($\alpha ,t>0$);
\item for fixed $t\in \mathbb{R}_+$, the function $\omega _{k}(t,x,B)$ is continuous in $x$.
\end{enumerate}
\par
Further, we establish an inequality of the Bernstein – Nikol’skii
type.
\begin{lemma}~Let $G(\lambda )$ be a nonnegative even function on
$\mathbb{R}$ that is nondecreasing on $\mathbb{R}_{+}$, let $x\in
\mathfrak{E}(B)$ and let $\sigma (x,B)\le \alpha $. Then
\begin{equation}
\left\Vert \Delta _{h}^{k}G(B)x\right\Vert \leq
h^{k}\alpha ^{k}G(\alpha )\left\Vert x\right\Vert ,\quad h>0,\, \,
k\in \mathbb{N}_{0}.\label{eq:1}
\end{equation}
\label{Lemma_1}\end{lemma}

\begin{proof} Since
\[
\sigma (x,B)\le \alpha\quad  \text{and}\quad \left|1-e^{i\lambda
h}\right|^{2k}=4^{2}\sin ^{2k}\frac{\lambda h}{2}\le \lambda
^{2k}h^{2k},\quad\lambda \in \mathbb{R},
\]
on the basis of operational calculus for the operator $B$ we get
\begin{equation}
\begin{gathered}
\left\Vert \Delta _{h}^{k}G(B)x\right\Vert ^{2}=\int _{-\alpha }^{\alpha }\left|(1-e^{i\lambda h})^{k}\right|^{2}G^{2}(\lambda )d(E_{\lambda }x,x)\leq  \\
\le h^{2k}\int _{-\alpha }^{\alpha }\lambda ^{2k}G^{2}(\lambda
)d(E_{\lambda }x,x)\leq h^{2k}\alpha ^{2k}G^{2}(\alpha )\left\Vert
x\right\Vert ^{2}.  \label{thm1_eq1}
\end{gathered}
\end{equation}
 \end{proof}

For $k=0$ Lemma \ref{Lemma_1} yields

\begin{equation}
\left\Vert G(B)x\right\Vert \le G(\alpha )\left\Vert x\right\Vert
.\label{lem1_nasl0}
\end{equation}

\begin{cor}
Under the conditions of Lemma \ref{Lemma_1} with respect to $x$ and
$\sigma(x, B)$, the following relation is true:
\[ \| \Delta _{h}^{k}x\|
\le h^{k}\cdot \alpha ^{k}\cdot \| x\| ,\quad h\geq 0.\]
\label{thm1_nasl2}
\end{cor}

\begin{proof}
For the proof of this statement, it suffices to take $G(\cdot)
\equiv 1,\, \lambda\in \mathbb{R}$, in Lemma \ref{Lemma_1}.
\end{proof}

If $\mathfrak{H}=L_{2}([0,2\pi ])$ and $(Bx)(t)=ix'(t)$
\[
\mathcal{D}(B)=\left\{ x(t)\, |\, x\in W_{2}^{1}([0,2\pi ]),\,
x(0)=x(2\pi )\right\},
\]
where $W_{2}^{1}([0,2\pi ])$ is a Sobolev space, then
$\mathfrak{E}(B)$ coincides with the set of all trigonometric
polynomials, $\sigma (x,B)$ is the degree of the polynomial $x$,
$\mathfrak{E}^{\alpha }(B)$ is the set of all trigonometric
polynomials whose degrees do not exceed $\alpha $;
$(U(h)x)(t)=\widetilde{x}(t+h)$, $\omega _{k}(t,x,B)$ is the $k$th
modulus of continuity of the function $x(t)$, and inequality
(\ref{eq:1}) for $G(\lambda )=\left|\lambda ^{m}\right|$ and $k=0$
turns into a Bernstein-type inequality in the space $L_{2}[0,2\pi ]$
\cite{Axiezer_Appr}~ (here $\widetilde{x}(t)$ is understood as the
$2\pi $-periodic extension of the function $x(t)$).

For an arbitrary $x\in \mathfrak{H}$ following
\cite{MGorb_BDVectors,MGorb_OperAppr}, we set
\[
\mathcal{E}_{r}(x,B)=\inf _{y\in \mathfrak{E}(B)\, :\, \sigma
(y,B)\leq r}\left\Vert x-y\right\Vert ,\quad r>0,
\]
i.e., $\mathcal{E}_{r}(x,B)$ is the best approximation of the
element $x$ by exponential-type entire vectors $y$ of the operator
$B$ for which $\sigma (y,B)\leq r$. For fixed $x$,~
$\mathcal{E}_{r}(x,B)$ does not increase and
$\mathcal{E}_{r}(x,B)\rightarrow 0,\, \, r\rightarrow \infty $.~ It
is clear that
\[
\mathcal{E}_{r}(x,B)=\left\Vert
x-E([-r,r])x\right\Vert =\left\Vert x-F([0,r])x\right\Vert ,
\]
where $F(\Delta )$ is the spectral measure of the operator
$|B|=\sqrt{B^{*}B}$.

\begin{theorem}
Suppose that $G(\lambda )$ satisfies the conditions of Lemma
\ref{Lemma_1}. Then, for any $x\in \mathcal{D}(G(B))$ the following
relation is true:
\begin{equation}
\forall \, k\in \mathbb{N}\quad
\mathcal{E}_{r}(x,B)\leq \frac{\sqrt{k+1}}{2^{k}G(r)}\omega
_{k}\left(\frac{\pi }{r},G(B)x,B\right),\quad
r>0.\label{eq:2}
\end{equation}
 \label{Th2}
\end{theorem}

\begin{proof}
Using the spectral representation for the operator $B$ and the
monotonicity of the function $G(\lambda )$, we obtain
\begin{gather*}
\omega _{k}^{2}(t,G(B)x,B)=\underset{0<\tau \le t}{\sup }\left\Vert (e^{i\tau B}-\mathbb{I})^{k}G(B)x\right\Vert ^{2}\geq \left\Vert (e^{itB}-\mathbb{I})^{k}G(B)x\right\Vert ^{2}=  \\
=\int _{-\infty }^{\infty }|e^{i\lambda t}-1|^{2k}G^{2}(\lambda )d(E_{\lambda }x,x)=2^{k}\int _{\mathbb{R}}(1-\cos \lambda t)^{k}G^{2}(\lambda )d(E_{\lambda }x,x)\geq   \\
\geq 2^{k}G^{2}(r)\int _{|\lambda |\geq r}(1-\cos \lambda
t)^{k}d(E_{\lambda }x,x).
\end{gather*}
 We fix $r>0$ and take $t:\  0\le t\le \frac{\pi }{r}$. Then
$\sin rt\ge 0$. We multiply both sides of the above inequality by
$\sin rt$ and integrate the result with respect to $t$ from $0$ to
$\frac{\pi }{r}$. Then
\begin{equation}
\begin{gathered}
\int _{0}^{\pi /r}\omega _{k}^{2}(t,G(B)x,B)\sin rt\, dt\ge 2^{k}G^{2}(r)\int _{0}^{\pi /r}\int _{|\lambda |\geq r}(1-\cos \lambda t)^{k}\sin rt\, d(E_{\lambda }x,x)dt=  \\
=2^{k}G^{2}(r)\int _{|\lambda |\ge r}\left(\int _{0}^{\pi /r}(1-\cos
\lambda t)^{k}\sin rt\, dt\right)d(E_{\lambda }x,x).
\label{thm2_eq1}
\end{gathered}
\end{equation}

Since the function $\omega _{k}^{2}(t,G(B)x,B)$ is monotonically
nondecreasing, we have
\begin{equation} \int _{0}^{\pi /r}\omega _{k}^{2}(t,G(B)x,B)\sin
rt\, dt\le \int _{0}^{\pi /r}\omega _{k}^{2}\left(\frac{\pi
}{r},G(B)x,B\right)\sin rt\, dt=\frac{2}{r}\omega
_{k}^{2}\left(\frac{\pi
}{r},G(B)x,B\right).\label{thm2_eq2}
\end{equation}

Using the inequality (see \cite{Stepanets_Serduyk})
\begin{equation}
\int _{0}^{\pi }(1-\cos \theta t)^{k}\sin t\, dt\ge
\frac{2^{k+1}}{k+1},\quad \theta \ge 1,\ k\in
\mathbb{N}\label{thm2_eq3}
\end{equation}
and relations (\ref{thm2_eq1}) and (\ref{thm2_eq2}), we get
\begin{equation}
\frac{2}{r}\omega _{k}^{2}\left(\frac{\pi }{r},G(B)x,B\right)\geq
2^{k}G^{2}(r)\int _{|\lambda |\ge
r}\left(\frac{1}{r}\frac{2^{k+1}}{k+1}\right)d(E_{\lambda
}x,x)=\frac{2^{2k+1}G^{2}(r)}{r(k+1)}\mathcal{E}_{r}^{2}(x,B),
\label{thm2_eq4}
\end{equation}
 which is equivalent to \eqref{eq:2}.
\end{proof}

For $G(\lambda )=|\lambda |^{m}$, $\lambda \in \mathbb{R}$, $m>0$
Theorem \ref{Th2} yields the following corollary:

\begin{cor}
Let $x\in \mathcal{D}(|B|^{m}),\ m>0$. Then, for any  $k\in
\mathbb{N}$
\begin{equation}
\mathcal{E}_{r}(x,B)\leq
\frac{\sqrt{k+1}}{2^{k}r^{m}}\omega _{k}\left(\frac{\pi
}{r},|B|^{m}x,B\right),\quad r>0.\label{eq:3}
\end{equation}
\end{cor}

For the case where $B$ is the operator of differentiation with
periodic boundary conditions in the space
$\mathfrak{H}=L_{2}([0,2\pi ])$, i.e., $(Bx)(t)=ix'(t)$ and
$\mathcal{D}(B)=\left\{ x(t)\, |\, x\in W_{2}^{1}([0,2\pi ]),\,
x(0)=x(2\pi )\right\} $, inequality (\ref{eq:3}) is presented in
\cite{Chernux} for $k=1$ and in \cite{Stepanets_Serduyk} for
arbitrary $k\in \mathbb{N}$.

We now formulate the inverse theorem in the case of approximation of
a vector $x$ by exponential-type entire vectors of the operator $B$.

\begin{theorem} Let $\omega (t)$ be a function of the type of a modulus of continuity for which the following conditions
are satisfied:
\begin{description}
\item [1)]$\omega (t)$ is continuous and nondecreasing for $t\in \mathbb{R}_{+}$;
\item [2)]$\omega (0)=0$;
\item [3)]$\exists c>0\, \, \forall t>0\quad \omega (2t)\leq c\, \omega (t)$;
\item [4)]${\displaystyle \int _{0}^{1}}\frac{\omega (t)}{t}dt<\infty $.
\end{description}
Also assume that the function $G(\lambda )$ is even, nonnegative,
and nondecreasing for $\lambda \geq 0$, and, furthermore,
$\underset{\lambda
>0}{\sup }\frac{G(2\lambda )}{G(\lambda )}<\infty $.

If, for $x\in \mathfrak{H}$, there exists $m>0$ such that
\begin{equation}
\mathcal{E}_{r}(x,B)<\frac{m}{G(r)}\omega
\left(\frac{1}{r}\right),\quad r>0,\label{thm3_eq2}
\end{equation}
 then $x\in \mathcal{D}(G(B))$ and, for every $k\in \mathbb{N}$,
there exists a constant $m_{k}>0$ such that
\begin{equation}
\omega
_{k}(t,G(B)x,B)\leq m_{k}\left[t^{k}\int \limits
_{t}^{1}\frac{\omega (\tau )}{\tau ^{k+1}}d\tau +\int \limits
_{0}^{t}\frac{\omega (\tau )}{\tau }d\tau \right],\quad 0<t\leq
\frac{1}{2}.\label{thm3_eq3}
\end{equation}
\label{Th3}
\end{theorem}

First, we prove the following statement:

\begin{lemma}

Suppose that the function $\omega (t)$ satisfies conditions 1),2),3)
of Theorem \ref{Th3}. If, for $x\in \mathfrak{H}$, there exists
$c>0$ such that
\begin{equation}
\mathcal{E}_{r}(x,B)<m\omega
\left(\frac{1}{r}\right),\quad
r>0\label{thm3_{l}em1_{e}q2}
\end{equation}
 then, for every $k\in \mathbb{N}$, there exists a constant $c_{k}>0$
such that
\begin{equation}
\omega _{k}(t,x,B)\leq c_{k}\cdot
t^{k}\int \limits _{t}^{1}\frac{\omega (\tau )}{\tau ^{k+1}}d\tau
,\quad 0<t\leq \frac{1}{2}.\label{thm3_{l}em1_{e}q3}
\end{equation}
 \label{thm3_lemma1}
\end{lemma}
\begin{proof} It follows from condition \eqref{thm3_{l}em1_{e}q2}
that there exists a sequence $\{u_{2^{i}}\}_{i=0}^{\infty }$ of
exponential-type entire vectors such that $\sigma (u_{2^{i}},B)\leq
2^{i}$ and
\begin{equation}
\left\Vert x-u_{2^{i}}\right\Vert \leq
m\cdot \omega
\left(\frac{1}{2^{i}}\right).\label{thm3_{l}em1_{e}q4}
\end{equation}

We take an arbitrary $h\in \left(0,\frac{1}{2}\right]$ and choose a
number $N$ so that $\frac{1}{2^{N+1}}<h\leq \frac{1}{2^{N}}$.
Inequality \eqref{thm3_{l}em1_{e}q4} yields
\begin{equation}
\Delta _{h}^{k}x=\Delta _{h}^{k}u_{1}+\sum _{j=1}^{N}\Delta
_{h}^{k}\big (u_{2^{j}}-u_{2^{j-1}}\big )+\Delta _{h}^{k}\big
(x-u_{2^{N}}\big )\label{thm3_{l}em1_{e}q5}
\end{equation}
одержуємо:
\begin{equation}
\begin{gathered}
 \left\Vert u_{2^{j}}-u_{2^{j-1}}\right\Vert \leq \left\Vert u_{2^{j}}-x\right\Vert +\left\Vert x-u_{2^{j-1}}\right\Vert \leq  \\
 \leq m\cdot \omega \left(\frac{1}{2^{j}}\right)+m\cdot \omega \left(\frac{1}{2^{j-1}}\right)\leq 2m\cdot \omega \left(\frac{1}{2^{j-1}}\right)\leq 2c\,m\cdot \omega \left(\frac{1}{2^{j}}\right).\label{thm3_{l}em1_{e}q7}
 \end{gathered}
\end{equation}

By virtue of the monotonicity of $\omega (t)$, we have
\begin{equation}
\begin{gathered}
2^{k}\int _{1/2^{j}}^{1/2^{j-1}}\frac{\omega (u)}{u^{k+1}}du\geq
2^{k}\omega \left(\frac{1}{2^{j}}\right)\int
_{1/2^{j}}^{1/2^{j-1}}\frac{1}{u^{k+1}}du=\frac{2^{kj}}{k}\omega
\left(\frac{1}{2^{j}}\right)(2^{k}-1)\geq 2^{kj}\omega
\left(\frac{1}{2^{j}}\right).  \label{thm3_{l}em1_{e}q8}
\end{gathered}
\end{equation}

Since $\sigma(u_{2^{j}}-u_{2^{j-1}}, B)\le 2^j$ and $\sigma(u_1,
B)\le 1$, according to Corollary \ref{thm1_nasl2} we get
\begin{align*}
\Vert \Delta _{h}^{k}u_{1}\Vert &\leq h^{k}\cdot \Vert u_{1}\Vert ,\label{thm3_lem1_eq10}\\
\Vert \Delta _{h}^{k}(u_{2^{j}}-u_{2^{j-1}})\Vert &\leq h^{k}\cdot
(2^{j})^{k}\Vert u_{2^{j}}-u_{2^{j-1}}\Vert .
\end{align*}

Relations \eqref{thm3_{l}em1_{e}q4}, \eqref{thm3_{l}em1_{e}q7} and
(\ref{thm3_{l}em1_{e}q8}) yield
\begin{equation*}
\left\Vert \Delta _{h}^{k}(u_{2^{j}}-u_{2^{j-1}})\right\Vert \leq 2c
m h^{k}\cdot 2^{kj}\omega \left(\frac{1}{2^{j}}\right)\leq
2^{k+1}cmh^{k}\int _{1/2^{j}}^{1/2^{j-1}}\frac{\omega
(u)}{u^{k+1}}du\label{thm3_lem1_eq9}
\end{equation*}
and
\[
\left\Vert \Delta
_{h}^{k}(x-u_{2^{N}})\right\Vert \leq (\left\Vert e^{ihB}\right\Vert
+1)^{k}\left\Vert x-u_{2^{N}}\right\Vert \leq 2^{k}\cdot \left\Vert
x-u_{2^{N}}\right\Vert \leq 2^{k}m\cdot\omega
\left(\frac{1}{2^{N}}\right).
\]

Using these inequalities, we obtain
\begin{multline*}
\left\Vert \Delta _{h}^{k}x\right\Vert = \bigg\|\Delta_h^ku_0+\sum_{j=1}^N\Delta_h^k(u_j-u_{j-1})+\Delta_h^k(x-u_N)\bigg\|\le\\
\leq h^{k}\left\Vert u_{0}\right\Vert +2^{k+1}cmh^{k}\sum _{j=1}^{N}\int _{1/2^{j}}^{1/2^{j-1}}\frac{\omega (u)}{u^{k+1}}du+2^{k}m\cdot \omega \left(\frac{1}{2^{N}}\right)\leq  \\
\leq h^{k}\left\Vert u_{0}\right\Vert +2^{k+1}cmh^{k}\int
_{1/2^{N}}^{1}\, \frac{\omega (u)}{u^{k+1}}du+2^{k}m\cdot \omega
(2h)\leq\\
\leq h^{k}\left\Vert u_{0}\right\Vert +2^{k+1}cmh^{k}\int _{h}^{1}\, \frac{\omega (u)}{u^{k+1}}du+2^{k}cm\cdot \omega (h) = \\
=h^{k}\left(\left\Vert u_{0}\right\Vert +2^{k+1}cm\int _{h}^{1}\, \frac{\omega (u)}{u^{k+1}}du+2^{k}cm\frac{k}{1-h^{k}}\, \int _{h}^{1}\, \frac{\omega (h)}{u^{k+1}}du\right)\leq  \\
\leq c_{k}\cdot h^{k}\int _{h}^{1}\frac{\omega (u)}{u^{k+1}}du,\quad
\quad \text{where}\ c_{k}=\frac{\left\Vert u_{0}\right\Vert }{\int
_{1/2}^{1}\, \frac{\omega (u)}{u^{k+1}}du}+2^{k+1}cm+
2^{k}cm\frac{k}{1-\frac{1}{2^{k}}}.\tag*{\hspace{-1em}\qedhere}
\end{multline*}
\end{proof}

\begin{rmk}

As follows from the proof, the lemma remains true under somewhat
weaker conditions than those formulated in the theorem, namely, it
is sufficient that, for an element $x\in \mathfrak{H}$, there exist
at least one sequence $\{u_{2^{j}}\}_{j=0}^{\infty }$, such that
\[
\sigma (u_{2^{j}},B)\leq 2^{j} \quad\text{and}\quad\forall j\in
\mathbb{N}\ \left\Vert x-u_{2^{j}}\right\Vert \leq m\cdot \omega
\left(\tfrac 1{2^{j}}\right).
\] \label{Lemma_remark}
\end{rmk}
\begin{proof}[Proof of Theorem] By virtue of \eqref{thm3_eq2}
there exists a sequence $\left\{ u_{2^{n}}\right\} _{n=1}^{\infty }$
such that $\sigma (u_{2^{n}})\leq 2^{n}$ and
\begin{equation}
\Vert
x-u_{2^{n}}\Vert \leq \frac{c}{G(2^{n})}\omega \left(\tfrac
1{2^{n}}\right),\quad n\in \mathbb{N}.\label{thm3_eq4}
\end{equation}

It follows from inequality (\ref{thm3_eq4}) and conditions 1), 2) of
the theorem that $\left\Vert x-u_{2^{n}}\right\Vert \rightarrow 0$
as $n\to \infty $, and, therefore, the vector $x$ can be represented
in the form
\[
x=u_{1}+\sum _{k=1}^{\infty }\big (u_{2^{k}}-u_{2^{k-1}}\big ).
\]

Since $\sigma (u_{2^{k}}-u_{2^{k-1}},B)\leq 2^{k}$, $k\in
\mathbb{N}$ taking (\ref{lem1_nasl0}) into account we obtain
\begin{gather*}
\left\Vert G(B)u_{2^{k}}-G(B)u_{2^{k-1}}\right\Vert \leq G(2^{k})\left\Vert u_{2^{k}}-u_{2^{k-1}}\right\Vert \leq G(2^{k})\left(\left\Vert x-u_{2^{k}}\right\Vert +\left\Vert x-u_{2^{k-1}}\right\Vert \right)\leq  \\
\leq G(2^{k})\left(\frac{m}{G(2^{k})}\omega \left(\frac{1}{2^{k}}\right)+\frac{m}{G(2^{k-1})}\omega \left(\frac{1}{2^{k-1}}\right)\right)\leq  \\
\leq \frac{2G(2^{k})\cdot m}{G(2^{k-1})}\omega
\left(\frac{1}{2^{k-1}}\right)\leq 2cc_{1}m\cdot \omega
\left(\frac{1}{2^{k}}\right)\leq\frac{2cc_1m}{\ln 2}\int
_{2^{-k}}^{2^{-k+1}}\frac{\omega (u)}{u}\, du,
\end{gather*}
where $c_1$ denotes $\underset{\lambda >0}{\sup }\frac{G(2\lambda
)}{G(\lambda )}$. Therefore, the series $\sum _{k=1}^{\infty }\big
(G(B)u_{2^{k}}-G(B)u_{2^{k-1}}\big )$ converges. The closedness of
the operator $G(B)$ implies that $x\in \mathcal{D}(G(B))$ and
\[
G(B)x=G(B)u_{1}+\sum _{k=1}^{\infty }\big
(G(B)u_{2^{k}}-G(B)u_{2^{k-1}}\big ).
\]
 This yields
\begin{gather*}
\left\Vert G(B)x-G(B)u_{2^{j}}\right\Vert \leq \sum _{k=j+1}^{\infty }\left\Vert G(B)u_{2^{k}}-G(B)u_{2^{k-1}}\right\Vert \leq 2cc_{1}m\sum _{k=j+1}^{\infty }\omega (2^{-k})\leq  \\
\leq 2cc_{1}m\int _{0}^{2^{-j}}\frac{\omega (u)}{u}\,
du=:\widetilde{c}\, \Omega (2^{-j}),\quad j\in \mathbb{N}
\end{gather*}
 where
 \[
 \tilde{c}:=2cc_{1}m\quad\text{and}\quad \Omega (t):=\int _{0}^{t}\frac{\omega (u)}{u}\, du
 \]
 It is easy to verify that the function $\Omega (t)$ possesses the following properties:

\begin{description}
\item [1)]$\Omega (t)$ is continuous and monotonically nondecreasing;
\item [2)]$\Omega (0)=0$;
\item [3)] for $t>0$, the following relation is true:
\[
\Omega (2t)=\int _{0}^{2t}\frac{\omega (u)}{u}\, du=\int
_{0}^{t}\frac{\omega (2u)}{u}\, du\leq c_{2}\int
_{0}^{t}\frac{\omega (u)}{u}\, du=c_{2}\Omega (t).
\]
\end{description}
Therefore, setting $\omega (t):=\Omega (t)$ in Lemma
\ref{thm3_lemma1} and taking Remark \ref{Lemma_remark} into account,
we get
\begin{gather*}
\omega _{k}\big (G(B)x,t,B\big )\leq c_{k}\cdot t^{k}\int _{t}^{1}\frac{\Omega (u)}{u^{k+1}}\, du=\frac{c_{k}\cdot t^{k}}{k}\left(\Omega (u)\left.\frac{1}{u^{k}}\right|_{1}^{t}+\int _{t}^{1}\frac{\omega (u)}{u^{k+1}}\, du\right)\leq   \\
\leq m_{k}\left(t^{k}\int _{t}^{1}\frac{\omega (u)}{u^{k+1}}\,
du+\int _{0}^{t}\frac{\omega (u)}{u}\,
du\right).\tag*{\hspace{-1em}\qedhere}
\end{gather*}
\end{proof}

Theorem \ref{Th3} shows that, in the case where $\omega
(t)=t^{\alpha },\, \, t\geq 0$, $\alpha
>0$ and $\mathcal{E}_{r}(x,B)=O\left(\frac{1}{r^{\alpha }}\right)$, one has
\[
\omega _{k}(t,x,B)=\left\{ \begin{array}{ll}
 O\left(t^{k}\right) & \texttt {при}\, \, k<\alpha \\
 O\left(t^{k}|\ln t|\right) & \texttt {при}\, \, k=\alpha \\
 O\left(t^{\alpha }\right) & \texttt {при}\, \, k>\alpha \quad .\end{array}\right.
 \]
 \par{}

2. Consider the equation
\begin{equation}
Ax=y,\label{eq:4}
\end{equation}
 where $A$ is a positive-definite self-adjoint operator with discrete spectrum, $y\in \mathfrak{H}$,
 $x\in \mathcal{D}(A)$ is the required
solution of Eq. (\ref{eq:4}). Let $\mathfrak{H}_{+}$ denote the
completion of the set $\mathcal{D}(A)$ with respect to the norm
$\left\Vert \cdot \right\Vert _{+}$, generated by the scalar product
\[
(x,y)_{+}=(Ax,y)\, .
\]
 Under the conditions imposed above on the operator $A$, Eq. (\ref{eq:4})
has a unique solution $x\in \mathcal{D}(A)$ and, according to the
Dirichlet principle \cite{Mixlin}, the determination of this
solution is equivalent to the determination of the vector $u\in
\mathcal{D}(A)$, on which the functional
\[
F(z)=(Az,z)-2Re(y,z),
\]
 defined on $\mathcal{D}(A)$ attains its minimum.

Let $\left\{ e_{k}\right\} _{k=1}^{\infty }$ be a complete linearly
independent system of vectors from $\mathcal{D}(A)$ (so-called
coordinate system), and let
\[
\mathcal{H}_{n}=\texttt {Л.О.}\left\{
e_{1},\cdots ,e_{n}\right\} .
\]
By $x_{n}$ we denote the vector on which $F(z)$ attains its minimum
on $\mathcal{H}_{n}$. The vector $x_{n}$ is called the Ritz
approximate solution of Eq. (\ref{eq:4}). As is known, independently
of the choice of a coordinate system, the sequence $x_{n}$ converges
to $x$ in the space $\mathfrak{H}_{+}$ (and, hence, in
$\mathfrak{H}$). The residual $R_{n}=\left\Vert Ax_{n}-y\right\Vert
$ does not always tend to zero in $\mathfrak{H}$. However, if the
coordinate system $\left\{ e_{k}\right\} _{k=1}^{\infty }$ is chosen
so that it forms an orthonormal proper basis of some
positive-definite self-adjoint operator $B$ related to $A$ in the
sense that $\mathcal{D}(A)=\mathcal{D}(B)$, then $R_{n}\rightarrow
0$ as $n\rightarrow \infty $ (see \cite{Mixlin}), and, therefore,
the quantities $r_{n}=\left\Vert x_{n}-x\right\Vert _{+}$ also tend
to zero as $n\rightarrow \infty $. However, the investigation of the
behavior of these quantities, which depend on the choice of $\left\{
e_{k}\right\} _{k=1}^{\infty }$ and on the right-hand side of Eq.
(\ref{eq:4}), at infinity turned out to be a rather difficult
problem and remains unsolved. Some particular results for operators
generated by boundary-value problems for ordinary differential
equations were obtained in numerous papers by many authors (see the
survey \cite{Lychka}). For the abstract case, some particular
situations were considered in \cite{Djushkariani}). In
\cite{MGorb_OperAppr}, direct and inverse theorems were established
for the first time under the condition that $x\in C^{\infty }(B)$
and estimates for the quantity $R_{n}$ were obtained in the case
where the smoothness of the vector  $x$ is finite, i.e., $x\in
\mathcal{D}(B^{k})$. Below, we completely characterize the quantity
$r_{n}$ for $x\in \mathcal{D}(B^{k})$.

In what follows, we assume that the following conditions are
satisfied:
\begin{description}
\item [$1^{0}$] The operator $A$ is self-adjoint and positive definite.
\item [$2^{0}$] The coordinate system in the Ritz method is an orthonormal basis of a positive-definite self-adjoint
operator $B$ with discrete simple spectrum ($Be_{k}=\lambda
_{k}e_{k}$) that is related to $A$.
\end{description}

Let $x_{n}$ denote the Ritz approximate solution of Eq. (\ref{eq:4})
with respect to the coordinate system $\left\{ e_{k}\right\}
_{k=1}^{\infty }$. We set
\[
\widetilde{x}_{n}={\displaystyle \sum _{k=1}^{n}(x,e_{k})e_{k}}.
\]
Since the operators $A$ and $B$ are positive definite and
self-adjoint and $\mathcal{D}(A)=\mathcal{D}(B)$, it follows from
the Heinz inequality \cite{Birman_Solomjak} that
$\mathcal{D}\left(A^{\alpha }\right)=\mathcal{D}\left(B^{\alpha
}\right)$ for any $\alpha \in (0,1)$, and, therefore, the operators
$B^{\frac{1}{2}}A^{-\frac{1}{2}}$ and
$A^{\frac{1}{2}}B^{-\frac{1}{2}}$ are defined and bounded on the
entire space $\mathfrak{H}$, and, for any $x\in \mathcal{D}(A)$, one
has
\begin{equation}
\mathbf{c}_{1}^{-1}|||x|||_{+}\leq
||x||_{+}\leq \mathbf{c}_{2}|||x|||_{+}\,,\label{eq:5}
\end{equation}
where $|||x|||_{+}=\left\Vert B^{1/2}x\right\Vert $,
$\mathbf{c}_{1}=\left\Vert B^{1/2}A^{-1/2}\right\Vert $ and
$\mathbf{c}_{2}=\left\Vert A^{1/2}B^{-1/2}\right\Vert $.

\begin{lemma}
For any $n\in \mathbb{N}$ and $x\in \mathcal{D}(B)$, the following
inequality is true:
\begin{equation}
|||x-\widetilde{x}_{n}|||_{+}\leq |||x-x_{n}|||_{+}\leq
\mathbf{c}_{3}|||x-\widetilde{x}_{n}|||_{+}\,,\label{eq:6}
\end{equation}
 where $\mathbf{c}_{3}=\left\Vert B^{1/2}A^{-1/2}\right\Vert \left\Vert A^{1/2}B^{-1/2}\right\Vert $.\label{Lemma_2}
\end{lemma}

\begin{proof}
Since
\[
B^{1/2}\big(\sum _{k=1}^{n}(x,e_{k})e_{k}\big)=\sum
_{k=1}^{n}\left(B^{1/2}x,e_{k}\right)e_{k},
\]
we have
\begin{gather*}
|||x-\widetilde{x}_{n}|||_{+}=\left\Vert B^{1/2}\left(x-\sum _{k=1}^{n}(x,e_{k})e_{k}\right)\right\Vert =\left\Vert B^{1/2}x-\sum _{k=1}^{n}\left(B^{1/2}x,e_{k}\right)e_{k}\right\Vert \leq  \\
\leq \left\Vert B^{1/2}x-B^{1/2}x_{n}\right\Vert =|||x-x_{n}|||_{+}
\end{gather*}

Taking into account that the Ritz approximation $x_{n}$ is the best
approximation of a vector $x$ in the norm $\left\Vert \cdot
\right\Vert _{+}$, we get
\begin{gather*}
|||x-x_{n}|||_{+}=\left\Vert B^{1/2}(x-x_{n})\right\Vert \leq \left\Vert B^{1/2}A^{-1/2}\right\Vert \left\Vert A^{1/2}(x-x_{n})\right\Vert =\mathbf{c}_{1}\left\Vert x-x_{n}\right\Vert _{+}\leq   \\
\leq \mathbf{c}_{1}\left\Vert x-\widetilde{x}_{n}\right\Vert
_{+}=\mathbf{c}_{1}\left\Vert
A^{1/2}\left(x-\widetilde{x}_{n}\right)\right\Vert \leq
\mathbf{c}_{1}\mathbf{c}_{2}\left\Vert
B^{1/2}\left(x-\widetilde{x}_{n}\right)\right\Vert
=\mathbf{c}_{3}|||x-\widetilde{x}_{n}|||_{+}\tag*{\hspace{-1em}\qedhere}
\end{gather*}
\end{proof}

Taking into account the relations
\[
\mathcal{E}_{\lambda
_{n}}(B^{1/2}x,B)=|||x-\widetilde{x}_{n}|||_{+}
\]
 and
 \[
\mathcal{E}_{\lambda _{n}}(B^{1/2}x,B)=\mathcal{E}_{\lambda
_{n}+\eta }(B^{1/2}x,B),\, \, 0<\eta <\lambda _{n+1}-\lambda _{n},
\]
 inequalities (\ref{eq:5}) and (\ref{eq:6}), and Theorem \ref{Th2}
 with
$G(\lambda )=|\lambda |^{\alpha -\frac{1}{2}}$, $\alpha \geq 1$, we
establish the following result:

\begin{theorem}
If $x\in \mathcal{D}(B^{\alpha }),\, \, \alpha \geq 1$, then the
following relation holds for every $\forall \, k\in \mathbb{N}$:
\[
||x-x_{n}||_{+}\leq \mathbf{c}_{0}\frac{\sqrt{k+1}}{2^{k}\, \lambda
_{n+1}^{\alpha -\frac{1}{2}}}\omega _{k}\left(\frac{\pi }{\lambda
_{n+1}},B^{\alpha }x,B\right),
\]
where $\mathbf{c}_{0}=\mathbf{c}_{2}\mathbf{c}_{3}$, and
$\mathbf{c}_{2}$ and $\mathbf{c}_{3}$ are the constants from
inequalities (\ref{eq:5}) and (\ref{eq:6}). \label{Th4}\end{theorem}

Since, for $x\in \mathcal{D}(B^{\alpha })$
\[
\omega
_{k}\left(\frac{\pi }{\lambda _{n+1}},B^{\alpha
}x,B\right)\rightarrow 0,\, \, n\rightarrow \infty ,
\]
 we conclude that, for $x\in \mathcal{D}(B^{\alpha })$
 \begin{equation}
\lim _{n\rightarrow \infty }\, \lambda _{n+1}^{\alpha
-\frac{1}{2}}||x-x_{n}||_{+}=0\label{eq:7}
\end{equation}

We now give examples of operators $A$ and $B$ for which equality
(\ref{eq:7}) for $\alpha >1$ does not yield the inclusion $x\in
\mathcal{D}(B^{\alpha })$. We set
\begin{gather*}
\mathfrak{H}=L_{2}([0,\pi ]),\quad A=B=-\frac{d^{2}}{dt^{2}},\quad
\mathcal{D}(A)=\mathcal{D}(B)=\left\{ x(t)\, |\, x\in
W_{2}^{2}([0,\pi ]),\, x(0)=x(\pi )=0\right\},\\
\lambda _{k}(B)=k^{2},\quad e_{k}=\sqrt{\frac{2}{\pi }}\sin kt,\quad
x=x(t)=\sqrt{\frac{2}{\pi }}{\displaystyle \sum _{k=2}^{\infty
}}x_{k}\sin kt,
\end{gather*}
where $x_{k}={\displaystyle \frac{1}{k^{2\alpha +\frac{1}{2}}\ln
^{\frac{1}{2}}k}}$, $k\in \mathbb{N}\backslash \{1\}$. The equality
\[
\sum _{k=2}^{\infty }\frac{k^{4\alpha }}{k^{4\alpha +1}\, \ln
k}=\sum _{k=2}^{\infty }\frac{1}{k\, \ln k}=\infty
\]
shows that $x\notin \mathcal{D}(B^{\alpha })$. However, since
\begin{gather*}
||x-x_{n}||_{+}^{2}=||x-\widetilde{x}_{n}||_{+}^{2}=\sum _{k=n+1}^{\infty }\frac{k^{2}}{k^{4\alpha +1}\, \ln k}\leq\\
\leq \frac{1}{\ln (n+1)}{\displaystyle \int \limits _{n}^{\infty
}}\frac{1}{t^{4\alpha -1}}dt=\frac{1}{(4\alpha -2)n^{4\alpha -2}\ln
(n+1)}
\end{gather*}
we have
\begin{align*}
\lim _{n\to \infty } & \lambda _{n}^{\alpha
-\frac{1}{2}}(B)||x-x_{n}||_{+}\leq \lim _{n\to \infty }n^{2\alpha
-1}\frac{1}{\sqrt{4\alpha -2}}\frac{1}{\sqrt{\ln (n+1)}n^{2\alpha
-1}}=0
\end{align*}

It follows from Theorem \ref{Th4}, inequality (\ref{eq:5}) and Lemma
\ref{Lemma_2} that the following statement is true:

\begin{theorem}
Suppose that $\omega (t)$ satisfies the conditions of Theorem
\ref{Th3}. If, for $x\in \mathcal{D}(B)$, $n\in \mathbb{N}$ and
$\alpha
>1$ one has
\[
\left\Vert x-x_{n}\right\Vert _{+}\leq
\frac{c}{\lambda _{n+1}^{\alpha -\frac{1}{2}}}\omega
\left(\frac{1}{\lambda _{n+1}}\right),
\]
 where $c\equiv \mathrm{c}onst$, then $x\in \mathcal{D}(B^{\alpha })$.
\label{Th5}
\end{theorem}

Note that, by virtue of inequality (\ref{eq:5}), $||\cdot ||_{+}$ in
Theorems \ref{Th4} and \ref{Th5} can be replaced by $|||\cdot
|||_{+}$.

The same theorem immediately yields the following corollary:

\begin{cor}
Suppose that the following inequality holds for $x\in
\mathcal{D}(B)$, $n\in \mathbb{N}$, $\alpha >1$ and $\varepsilon >0$
\[
\left\Vert x-x_{n}\right\Vert _{+}\leq
\frac{c}{\lambda _{n+1}^{\alpha +\varepsilon -\frac{1}{2}}}\, \, .
\]
 Then $x\in \mathcal{D}(B^{\alpha })$.
\end{cor}

\begin{rmk} If, as the Ritz approximate solution of (\ref{eq:4}), one takes the vector $x_{n}$
on which the functional $F(z)$ attains its minimum on
$\mathfrak{H}_{n}=\mathfrak{H}_{\lambda _{1}}\bigoplus
\mathfrak{H}_{\lambda _{2}}\bigoplus \cdots \bigoplus
\mathfrak{H}_{\lambda _{n}}$, where $\mathfrak{H}_{\lambda _{j}}$ is
the eigensubspace of the operator $B$ corresponding to the
eigenvalue $\lambda _{j}$, then, under assumption $2^{0}$ one can
omit the condition of the simplicity of the spectrum.
\end{rmk}

3. We set $\mathfrak{H}=L_{2}(0,\pi )$, $\mathcal{D}(A)=\left\{ x\in
W_{2}^{2}[0,\pi ],\, \, x'(0)=x'(\pi )=0\right\} $ and
\[
(Ax)(t)=-x''(t)+q(t)x(t),\quad q(t)>0,\, \, q\in C([0,\pi ]).
\]

We define an operator $B$ as follows:
\[
\mathcal{D}(B)=\mathcal{D}(A),\quad Bx=-x''+x.
\]

The operators $A$ and $B$ are self-adjoint and positive definite in
$L_{2}(0,\pi )$. The spectrum of $B$ consists of the eigenvalues
$\lambda _{k}(B)=k^{2}+1$, $k\in \mathbb{N}_{0}$, corresponding to
the eigenfunctions $\sqrt{\frac{2}{\pi }}\cos \left(kt\right)$,
which form an orthonormal basis in the space $L_{2}(0,\pi )$.

Let $k\in \mathbb{N}$ and $g(t)\in C^{2k}[0,2\pi ]$. It is easy to
verify that $\mathcal{D}(A^{k+1})=\mathcal{D}(B^{k+1})$ if and only
if $g^{2j+1}(0)=g^{2j+1}(\pi )=0$, $j=0,\cdots ,k$. If $y(t)\in
C^{2(k-1)}[0,2\pi ]$ and $y^{2j+1}(0)=y^{2j+1}(\pi )=0$, then
$y(t)\in \mathcal{D}(A^{k})$. Therefore, the solution of the problem
\begin{gather}
 -x''(t)+g(t)x(t)=y(t)\label{Zad_Koshi1}\\
 x'(0)=x'(\pi )=0\label{Zad_Koshi2}
\end{gather}
belongs to the set $\mathcal{D}(A^{k+1})=\mathcal{D}(B^{k+1})$ and
relation (\ref{eq:7}) directly yields the following statement:

\begin{theorem}
If $g(t)\in C^{2k}[0,\pi ]$, $g^{(2j+1)}(0)=g^{(2j+1)}(\pi )=0$,
$j=0,\cdots ,k$, and $y(t)\in C^{2(k-1)}[0,2\pi ]$,
$y^{(2j+1)}(0)=y^{(2j+1)}(\pi )=0$, $j=0,\cdots ,k-1$, then the Ritz
approximate solution of problem
(\ref{Zad_Koshi1})-(\ref{Zad_Koshi2}) satisfies the relation
\[
\left\Vert x_{n}-x\right\Vert _{W_{2}^{2}[0,\pi
]}=o\left(\frac{1}{n^{2k+1}}\right).
\]
\end{theorem}

 \end{document}